\documentclass[11pt,a4paper]{article}
\usepackage{dsfont}
\usepackage{parskip}
\usepackage{amssymb}
\begin{document}
\author{Anton D. Baranov, Ilgiz R. Kayumov, Semen R. Nasyrov}
\title{On Bloch seminorm of  finite Blaschke products in the unit disk}
\date{}
\makeatletter
\renewcommand{\@makefnmark}{}
\makeatother
\newtheorem{lemma}{Lemma}
\newtheorem{remark}{Remark}
\newtheorem{corollary}{Corollary}
\newtheorem{theorem}{Theorem}
\newtheorem{proposition}{Proposition}
\newtheorem{df}{Definition}

\maketitle
\begin{footnote}{Keywords: {Blaschke product,  Riemann surface, Bloch seminorm, bounded analytic function.}

AMS classification numbers: 30J10.}
\end{footnote}
\makeatletter
\renewcommand{\@makefnmark}{}

\makeatother

\maketitle
\begin{abstract}
We prove that, for any finite Blaschke product  $w=B(z)$ in the unit disk,
the corresponding Riemann surface over the $w$--plane contains a one-sheeted disk of the radius $0.5$.
 Moreover, it contains a unit one-sheeted disk with a radial slit.
 We apply this result to obtain a universal sharp lower estimate
 of the Bloch seminorm for finite Blaschke products.
\end{abstract}

\section{Introduction}

Let $\mathbb{D}$ be the unit disk $\{|z|<1\}$ in the complex plane
$\mathbb{C}$. We will consider finite Blaschke products which can be represented in the form
\begin{equation}\label{bl}
B(z)=\lambda \prod_{j=1}^n\frac{z-z_j}{1-\overline{z_j}z}, \qquad |z_j|<1,
\end{equation}
where $|\lambda|=1$; further, for simplicity of presentation of the results, we will assume that $\lambda=1$.
Finite Blaschke products provide very basic examples of bounded analytic functions
which have many important and nontrivial properties \cite{Garnett}, \cite{Gar},
 \cite{Mash}, \cite{KraRot}, \cite{Bog}.
Also, they have numerous applications in complex dynamics \cite{Mc}.

Recall that for a function $f$ analytic in the unit disc
$\mathbb{D}$ its Bloch seminorm is given by
$$
\|f\|_{\mathbb{B}}:=\sup_{z\in {\mathbb{D}}}|f'(z)|(1-|z|^2).
$$The class $\mathbb B$ of analytic functions with bounded Bloch seminorm is called the Bloch space. It is a Banach space if endowed with the norm $|f(0)| + \| f \|_{\mathbb B} $. 
The Bloch seminorm is closely related to the inner (conformal) radius of a domain $\Omega$ at a point $w_0 \in \Omega$ \cite{PolSeg}:
$$
r(w_0):=\frac{1-|g(w_0)|^2}{|g'(w_0)|}
$$ where $g$ conformally maps $\Omega$ onto $\mathbb{D}$; if additionally  we consider the mapping $g$ with  $g(w_0)=0$, then $r(w_0):=\frac{1}{|g'(w_0)|}$. One can easily check that $r(f(z_0))=(1-|z_0|^2)|f'(z_0)|$ where $f=g^{-1}$. The value $\sup_{w_0 \in \Omega}r(w_0)$ is called the maximal conformal radius of $\Omega$.

We also recall that if $h$  maps conformally $\Omega$ onto the
upper half-plane $\mathbb{H}$, then
\begin{equation}\label{upper}
r(w_0)=\frac{2\Im h(w_0)}{|h'(w_0)|}\,.
\end{equation}


By the Schwarz--Pick inequality \cite{AvkWir}, for any $f\in H^\infty$ (the
space of all functions bounded and analytic in $\mathbb{D}$) we
have $\|f\|_{\mathbb{B}} \le \|f\|_\infty$ and, therefore,
$\|B\|_{\mathbb{B}} \le 1$ for any (even infinite) Blaschke
product. On the other hand, if we consider a M\"obius
transformation, i.e. a Blaschke product of degree one
$$B_a(z)=\frac{z+a}{1+\overline{a}z}
$$
where $|a|<1$, then $\|B_a\|_{\mathbb{B}} =1$.
 For each $n$ the inequality $\|B\|_{\mathbb{B}} \le 1$ cannot be improved
 because the  Blaschke  product of degree one can be approximated locally uniformly  by  Blaschke  products of degree $n$. We want to investigate the following
\medskip
\\
{\bf Problem.} {\it Is there a universal constant  $c>0$ such
that for any finite Blaschke product $B$ the inequality
} $ \|B\|_{\mathbb{B}} \ge  c$   {\it holds?}
\medskip

From the results of Aleksandrov, Anderson  and Nicolau \cite{Alex}
it follows that, in the case of infinite Blaschke products, the
answer is negative, i.e. $c=0$. Hence, one may expect that the
value $$\inf\{\|B\|_{\mathbb{B}}: B \mbox{ is a Blaschke product
of degree } n\}$$ goes to zero as $n \to \infty$. Surprisingly, it
turned out that our problem has a positive answer. In addition it should be noted that finite Blaschke products belong to the little Bloch space
$$
\mathcal{B}_0=\{f \in \mathbb{B}: \lim_{|z| \to 1} |f'(z)|(1-|z|)=0\}.
$$

\begin{theorem}\label{Lower} Suppose that $B$ is a finite Blaschke product. Then
\begin{equation}\label{068}
\|B\|_{\mathbb{B}} \ge r_0=0.695356\ldots 
\end{equation}
Conversely, for any $k >2/e=0.73575\ldots$ there exists a finite
Blaschke product for which $\|B\|_{\mathbb{B}} \leq k$.
\end{theorem}

Theorem \ref{Lower} is closely related to classical covering
theorems in the geometric function theory. The classical Koebe
$1/4$--theorem (see, e.g. \cite[ch.~II, \S~4]{Goluzin}) implies
that, for every univalent analytic function $f: \mathbb{D}
\to \mathbb{C}$, the image $f(\mathbb{D})$ contains the disk with
center $f(0)$ and radius $|f'(0)|/4$. Now assume that $f$ is an
arbitrary analytic function $f: \mathbb{D} \to \mathbb{C}$. We
define $B_f$ to be the radius of the largest disk that is the
biholomorphic image of a subset of a unit disk. Bloch \cite{Blo}
proved that $B_f \ge |f'(0)|/72$. Ahlfors \cite{Ahl}, as an
application of his method of "ultrahyperbolic metrics",
established that $B_f \ge |f'(0)|\sqrt{3}/4$. Bonk \cite{Bonk}
slightly improved this estimate (see also \cite{Chen}).

The main goal of our paper is to establish an analog of Bloch's theorem for $f \circ B$ where is $B$ is a Blaschke product.
More precisely, we are going to investigate the following question: Let $f$ be an
analytic in the disk $\mathbb{D}$ function and $f'(0)=1$. Is it
true that there exists a universal constant $c$ such that  for all
Blaschke products $B$ there is a one-sheeted disk of the radius
$c$ that lies in the Riemann surface of the function inverse to
$f\circ B$. As mentioned above, the answer
is negative for infinite Blaschke products.

It turns out that such bound exists for finite Blaschke products.
To prove this, we need a statement on Riemann surfaces generated
by Blaschke products.

In this paper, by Riemann surface $R(f)$ we will mean a covering
$f:D\to G$, where $D$ and $G$ are some abstract Riemann surfaces
or even domains in the complex plane and $f$ is a holomorphic
function mapping $D$ onto $G$ (projection). In general, $R(f)$ is
either ramified or unramified and need not be unlimited, i.e. we
do not require that it has the curve lifting property \cite[p.~25,
Definition~4.13]{Forster}. 

{\bf Definition 1.} \it Two Riemann surfaces $f_k:D_k\to G$,
$k=1$, $2$, are called equivalent if there exists a biholomorphic
mapping $h$ of $D_1$ onto $D_2$ such that $f_2=h\circ f_1$.
\rm

 As a
rule, equivalent Riemann surfaces are not distinguished.

{\bf Definition 2.} \it  We will
say that a Riemann surface $f_2:D_2\to G_2$ contains a Riemann
surface $f_1:D_1\to G_1$ if there exists an injective holomorphic
function $h:D_1\to D_2$ such that $f_2=h\circ f_1$. 
\rm

If $D$ is a
planar domain and $id_D:D\to D$ is the identity mapping in $D$,
then we will identify the corresponding Riemann surface with the
domain $D$.

Every finite Blaschke product $B(z)$ generates the Riemann surface
$B:\mathbb{D}\to \mathbb{D}$ which is an $n$-sheeted covering; we
will denote it by $R(B)$. Theorem~\ref{Main} below describes some
subdomains contained in $R(B)$. To formulate the theorem, we need
to introduce some notation.

Let $0\le a<1$. Denote by $\mathbb{D}_a$ the unit disk
$\mathbb{D}$ with the slit along the segment $[a,1]$ of the real
axis.  Consider also the domain $G_a$ which is the union of the
half-strip
\begin{equation}\label{str}
S:=\{ \Re w<0,\ 0<\Im w<2 \pi \}
\end{equation}
and the rectangle
\begin{equation}\label{pia}
\Pi_a:=\{\log a<\Re w<0, \ 0<\Im w<3 \pi\},
\end{equation}
Let $G^*_a$ be the domain which is obtained by the reflection of
$G_a$ with respect to the real axis. Denote by $g:\mathbb{C}\to
\mathbb{C}$ the exponential mapping $g(w)=e^w$. Now consider the
two-sheeted Riemann surface $\mathbb{L}_a$ which is defined by the
non-ramified covering $g|_{G_a}:G_a\to g(G_a)$ and its reflection
$\mathbb{L}^*_a$ with respect to the real axis specified by the
covering $g|_{G^*_a}:G^*_a\to g(G^*_a)$. The surfaces
$\mathbb{L}_a$ and  $\mathbb{L}^*_a$ are obtained by gluing of the
slit unit disk $\mathbb{S}$ and a half of the annulus $a<|z|<1$
along the common boundary segment $[a,1]$.

Now consider the equation
$$
s\,\frac{1 - a(s)^2}{(s-a(s))^2}\,\frac{s-1}{s+1}\,=0.175,\quad
\mbox{\rm where} \quad a(s):=s\,\frac{2s-(s^2-1)}{2s+(s^2-1)}\,,
$$
on the segment $[1,1+\sqrt{2}]$. From Lemma~\ref{lemma1} below it
follows that it has a unique solution $s_0=2.379796\ldots$, and we
define $a=a(s_0)=0.024286\ldots$.


\begin{theorem} \label{Main}
Let $R(B)$ be the Riemann surface of a finite Blaschke product $B$
and $a=0.024286\ldots$. Then $R(B)$ contains, up to a rotation,
either the slit disk $\mathbb{D}_a$, or $\mathbb{L}_a$, or
$\mathbb{L}_a^*$.
\end{theorem}

From Theorem~\ref{Main} we obtain

\begin{corollary} \label{Main1} Let $B$ be a finite Blaschke product defined in the unit disk. Then $R(B)$ contains a disk of radius $1/2$. 
\end{corollary}


Theorem~$\ref{Main}$ allows us to investigate some properties of
Blaschke products in the Bloch space.
In fact,  Theorem~\ref{Lower} is a consequence of Theorem
\ref{Main}.

In the connection with Theorem~\ref{Lower}, we should note a
recent result by Dubinin~\cite[Theorem 1.1]{dub},  giving a sharp
upper bound of $|B'(z)|(1-|z|^2)$ for finite Blaschke products (\ref{bl}) provided
that the critical  values lie in a given disk.
\medskip

Theorem~\ref{Lower} allows
us to obtain similar results for functions generalizing Blaschke
products.

\begin{theorem} \label{General} Let $f$ be a holomorphic in $\mathbb{D}$ function with $f'(0)=1$ and $B$ be a finite Blaschke product.
Then the Bloch seminorm of $g=f\circ B$ satisfies the inequality
$$
\|g\|_{\mathbb{B}} \ge \sqrt{3}\,r_0/4  =0.301098\ldots%
$$
where $r_0$ is given in (\ref{068}). Moreover, if additionally $f$
is a convex univalent function, then
$$
\|g\|_{\mathbb{B}} \ge \pi r_0/4
=0.546131 \ldots
$$\end{theorem}

The proof of the lower bound in Theorem \ref{Lower} is not
constructive. We complement it with a simple observation how to
find a point $z$ where the quantity $|B'(z)|(1-|z|^2)$ admits a
universal (but smaller than in Theorem \ref{Lower}) lower bound.
This result also applies to some classes of infinite Blaschke
products.

\begin{theorem}\label{infi} Let  $B$ be a Blaschke product with zeros $\{z_j\}_{j \ge 1}$
and assume that there exists a point $\zeta\in \partial\mathbb{D}$ such
that ${\rm dist}\, (\zeta, \{z_j\}_{j \ge 1})>0$ and, moreover,
\begin{equation}
\label{onec}
|B'(\zeta)|\, {\rm dist}\,(\zeta, \{z_j\}_{j \ge 1}) \ge\delta>0
\end{equation}
for some $0< \delta \le 1$.
Then for $z_0=\big(1-\frac{\delta}{8|B'(\zeta)|}\big) \zeta$
we have
$$
|B'(z_0)|(1-|z_0|^2) \ge 0.07\delta.
$$
\end{theorem}
\medskip

If $B$ is a finite Blaschke product we can take $\zeta$ to be the
point where $|B'|$ attains its maximum in the closed disk
$\overline{\mathbb{D}}$. Then it is clear (see formula (\ref{pro})
below) that $|B'(\zeta)| \ge \frac{1+|z_j|}{1-|z_j|}$ and so
(\ref{onec}) is satisfied with constant $\delta = 1$.

Condition (\ref{onec}) appears, e.g., in the study of the so-called ``one-component'' Blaschke products
for which the level set $\{z\in\mathbb{D}: \ |B(z)| < \varepsilon \} $ is connected for some $ \varepsilon \in (0,1) $ (see \cite{b09, bd09}).
\bigskip

\section{Proofs of the main results}

To prove Theorem \ref{Main} we need to establish two lemmas.
\medskip

\begin{lemma}\label{lemma1} Let $0\le a<1$ and $\mathbb{D}_a$ be
the unit disk with the slit along the segment $[a,1]$ of the real
axis. Then the  maximal value of the conformal radius of
$\mathbb{D}_a$ is equal to the maximum of the function
$$
g(x)=\frac{ 4(1 -
    a^2)x(\sqrt{x^2 + 1} - x)^2 }{\sqrt{x^2 + 1} (1 - a
(\sqrt{x^2 + 1} - x)^2)^2},\quad x>0.
$$
\textit{The maximum is attained at the point
$x_0=(s_0-1)/(2\sqrt{s_0})$ where $s_0$ is the unique positive
root of the equation }
\begin{equation}\label{cub}
s^3-(2-a)s^2+(2a-1)s-a=0,\quad 1< s\le 1+\sqrt{2}.
\end{equation}
Moreover, 
$$
g(x_0)=4s_0\,\frac{1 -
    a^2}{(s_0-a)^2}\,\frac{s_0-1}{s_0+1}\quad \mbox{\rm and} \quad a=s_0\,\frac{2s_0-(s_0^2-1)}{2s_0+(s_0^2-1)}\,.
$$
\end{lemma}

{\bf Proof.} The conformal mapping of the lower half-plane onto
$\mathbb{D}_a$ has the form
$$
z=h(w)=\frac{(w-\sqrt{w^2-1})^2+a}{1+a(w-\sqrt{w^2-1})^2},
$$
consequently, the conformal radius of $\mathbb{D}_a$ at the point
$h(-ix)$, $x>0$, equals $2 x |h'(-ix)|=g(x)$.

The function $g$ is strictly positive for $x>0$ and
$$
\frac{g'(x)}{g(x)}\,=\,\frac{1}{x(x^2+1)}-\frac{2}{\sqrt{x^2+1}}\,\frac{1+a(\sqrt{x^2
+ 1} - x)^2}{1-a(\sqrt{x^2 + 1} - x)^2}\,,
$$
therefore, $g'(x)=0$ if and only if
\begin{equation}\label{eqq}
\frac{1-a(\sqrt{x^2 + 1} - x)^2}{1+a(\sqrt{x^2 + 1} -
x)^2}\,=2x\sqrt{x^2+1}.
\end{equation}
If $t=\log(\sqrt{x^2+1}+x)$, then $x=\sinh t$ and the equation
(\ref{eqq}) has the form
$$
\frac{e^{2t}-a}{e^{2t}+a}\,=\sinh 2t.
$$
It is equivalent to the cubic equation with respect to $s=e^{2t}$:
$$\psi(s)=0, \quad \mbox{\rm where} \quad \psi(s)=s^3-(2-a)s^2+(2a-1)s-a.
$$

Simple analysis shows that for real $s$ this equation has a unique
root $s_0$ satisfying the inequality $ 1< s_0\le 1+\sqrt{2}$.
Actually, $\psi(1)=2(a-1)<0$ and for $s=1+\sqrt{2}$ we have
$$
\psi(s)=s^3-2s^2-s+a(s^2+2s-1)=a(s^2+2s-1)\ge 0,
$$ therefore, there is at least one real root of
the cubic equation on the interval $(1,1+\sqrt{2}]$. Assume that
there is another real root $s_1$ on $(1,1+\sqrt{2}]$. By Vieta's
formulas, the product of the roots of the equation $\psi(s)=0$ is
equal $-a<0$, consequently, we have a third root $s_2$ of the
cubic equation which is negative. Then
$\psi(s)=(s-s_0)(s-s_1)(s-s_2)$ and, therefore, $\psi(1)>0$ but
this contradicts to the fact that $\psi(1)=2(a-1)<0$.

Thus, we have a unique root $s_0$ of (\ref{cub}). Then
$x_0=\sinh\frac{\log s_0}{2}=\frac{s_0-1}{2\sqrt{s_0}}$ is the
point of maximum of the function $g$. Lemma~\ref{lemma1} is
proved.
\medskip

\begin{remark}\label{rem1} {\rm From the proof of Lemma~\ref{lemma1} it follows that $g(x_0)$
can be expressed via $s_0$, i.e. $g(x_0)=F(s_0)$ for some
increasing function $F$. Calculations give that if $F(s_0) = 0.7$,
then $s_0=F^{-1}(0.7)=2.379796\ldots$,
\begin{equation}\label{s0}
 a=0.024286\ldots\ \ \mbox{\rm
and}\ \ \log a=-3.7178547\ldots.
\end{equation}}
\end{remark}

\begin{lemma}\label{lemma2}  Let $a$ be given by (\ref{s0}). The maximum of the conformal radius of the surface
$\mathbb{L}_a$ is not less than $ r_0:=0.695356\ldots$
\end{lemma}

{\bf Proof.} The conformal mapping of the upper half-plane
$\mathbb{H}$ onto $\mathbb{L}_a$ has the form $w=e^{f(z)}$. Here
$$
f(z)=-C\int_{-1}^z\frac{\sqrt{t-d}\,dt}{\sqrt{(t-c)(t^2-1)}}\,,\quad
1<c<d<+\infty,
$$
maps $\mathbb{H}$ onto $G_a$ which is the union of the strip $S$
and the rectangle $\Pi_a$, defined by (\ref{str}) and (\ref{pia}),
supplemented with their common boundary arc. The branches of the
square roots are fixed such that the integrand takes positive
values for real $t>d$, and the constant $C>0$.

Investigation of the behavior of $f(z)$ near $z=\infty$ gives
$$
f(z)=-C\log z+O(1)
$$
where $\log z$ takes positive values for real $z>0$. The width of
$S$ is equal to $2\pi$, therefore, we conclude that $C=2$. Since
the height and the width of the rectangle $\Pi_a$ equal $3\pi$ and
$-\log
a=3.7178547\ldots$ 
(see (\ref{s0})), we obtain
$$
\int_{-1}^1|f'(t)|dt=3\pi,\quad \int_{1}^c|f'(t)|dt=-\log a,
$$
therefore,
$$
\int_{-1}^1\frac{\sqrt{d-t}\,dt}{\sqrt{(c-t)(1-t^2)}}=3\pi/2,
$$
$$
\int_{1}^c\frac{\sqrt{d-t}\,dt}{\sqrt{(c-t)(t^2-1)}}=-\log
a/2=1.858927\ldots 
$$
Solving the system of equations with respect to $c$ and $d$ we
find 
$$
c=1.098259\ldots, 
\quad d=1.766556\ldots
$$

The conformal radius of $\mathbb{L}_a$ at the point $e^{f(z)}$ is
equal to
$$
2\Im z\left|\left(e^{f(z)}\right)'\right|=4\Im z\,
\left|e^{f(z)}\right| \frac{\sqrt{|z-d|}}{\sqrt{|(z-c)(z^2-1)|}}.
$$
At the point $z = -0.0205 + 0.3659i$ the conformal radius equals
$r_0=0.695356\ldots$ Therefore, Lemma~\ref{lemma2} is
proved.\medskip

\noindent {\bf Proof of Theorem \ref{Main}.}  Consider the surface
$R(B)$ for a finite Blaschke product. It is an $n$-sheeted
unlimited ramified covering of the unit disk $\mathbb{D}$. In the
case $n\le 2$ the statement of the theorem is evident, therefore,
we can assume that $n>2$.\medskip

I) First we consider the case where all branch points of $R(B)$
are simple and their projections $r_ke^{i\phi_k}$ ($\phi_k\in
[0,2\pi)$) on $\mathbb{D}$ are such that no two of the points lie
on the same radius of $\mathbb{D}$. Then $R$ can be glued from $n$
disks slit along the segments $T_k$ of the form $re^{i\phi_k}$,
$0<r_k<r<1$. We will call them sheets of the Riemann surface
$R(B)$ and denote them by $S_1,\ldots,S_n$. For every segment
$T_k$ there are exactly two sheets slit along it. From the
Riemann--Hurwitz formula for bordered surfaces (see, e.g.
\cite{nas}) it follows that the number of segments (and branch
points of $R(B)$) equals $(n-1)$. Since every segment corresponds
to a couple of slits, the number of slits on all the sheets equals
$2(n-1)$.

Now we will prove that there exists a sheet $S_k$ that contains
more than one slit and every sheet that is glued to it, with
possibly one exception, contains a unique slit. To prove this, we
associate with the surface $R(B)$, glued from the sheets $S_j$, a
connected graph $\Gamma$, the vertices of which are sheets. A
vertex $S_j$ is connected with $S_l$, if the sheets $S_j$ and
$S_l$ are glued with each other, with the help of the slit along
some segment $T_m$. It is easy to see that, because of
simply-connectedness  of $R(B)$, the graph is a tree. Now we
consider any (oriented) edge path of $\Gamma$ with the maximal
possible number of edges. Then the second vertex $S_k$ in this
path is the required sheet (and, if it exists, the exceptional
sheet with more than one slit is the third vertex in this path).

By renumbering the sheets and the points $r_j e^{i\phi_j}$, we can
achieve that $S_1$ is the sheet such that the sheets
$S_2,\ldots,S_m$, $m\le n-1$, are attached to it and each of them
has a unique slit along the segment connecting the
points $r_j e^{i\phi_j}$ and $e^{i\phi_j}$. 
We can assume that $0<\phi_1<\phi_2<\ldots<\phi_m<2\pi$. We can
glue from the sheets $S_1,\ldots,S_m$ a Riemann surface
$R_1\subset R(B)$ which has $m$ sheets, $m-1$ branch points and,
possibly, one more slit on some sheet $S_j$, $2\le j\le m$. We
will consider the case where there is such slit; in case of its
absence, we can always cut the surface along some radial segment.
Then we extend the slit so that it coincide with some radius.
Thus, the boundary of $R_1$ consists of the $m$ times traversed
unit circle and the slit.

If there is $k$, $1\le k\le m$, such that $r_k>a$, where $a$ is
given by (\ref{s0}), then the sheet $S_k$ contains a domain which
is a rotation of $\mathbb{D}_a$ by angle $\phi_k$ and the theorem
is proved. Therefore, we can assume that the projections of all
branch points of $R(B)$ are located at a distance from zero less
than $a$. Then $R_1$ contains a subsurface $R_2$ which is an
$m$-sheeted non-ramified covering of the annulus $a<|z|<1$, cut
along a radial segment. Without loss of generality we can assume
that the projection of the segment is on the positive part of the
real axis. Then, under an appropriate choice of the branch of the
logarithm, the function $w=\log z$ maps $R_2$ onto the rectangle
$\Pi=\{\log a<\Re w<0$, $0<\Im w<2 \pi m\}$. The union of $R_2$
and every  $S_k$, $1\le k\le m$, is a part of $R(B)$, containing
either $\mathbb{L}_a$ or $\mathbb{L}_a^*$, turned at an
angle~$\phi_k$. (We recall that we understand inclusion of one
Riemann surface into another one in the sense of Definition 2
given in Introduction).\medskip

II) Now let $B $ be an arbitrary Blaschke product. It can be
approximated by a sequence of Blaschke products $B_n$ satisfying
the requirements considered in~I) and the convergence is uniform
in the closed unit disk. Because of I), every $R(B_n)$ contains
either the slit disk $\mathbb{D}_a$, or $\mathbb{L}_a$, or
$\mathbb{L}_a^*$ rotated by some angle $\theta_n\in [0,2\pi]$.
Without loss of generality we can assume that
$\theta_n\to\theta_0$ as $n\to\infty$ and every $R(B_n)$ contains
a rotation of one of the indicated domains, $\mathbb{D}_a$,
$\mathbb{L}_a$, or $\mathbb{L}_a^*$. The sequence of  $R(B_n)$
converges to $R(B)$  as to a kernel in the sense of Caratheodory
(the kernel convergence of planar domain is described, e.g. in
\cite[\S~3.1, p.~77]{duren}, about the kernel convergence of
multi-sheeted Riemann surfaces see, e.g. \cite{nas1} and the
bibliography therein). Let for definiteness, $R(B_n)$ contains
$\mathbb{L}_a$ rotated by the angle $\theta_n$. Then $R(B)$
contains $\mathbb{L}_a$ rotated by the angle $\theta_0$.

Theorem~\ref{Main} is proved.
\medskip
\\
{\bf Proof of Theorem \ref{Lower}.} The norm $\|B\|_{\mathbb{B}}$
is the maximum of the conformal radius of the surface $R(B)$. By
Theorem~\ref{Main}, the surface $R(B)$ contains a rotation of
either $\mathbb{D}_a$, or $\mathbb{L}_a$, or $\mathbb{L}_a^*$. We
note that by the Lindel\"of principle (see \cite{Goluzin}, p. 339),
the conformal radius
increases under enlargement of domain (Riemann surface). Since the
maximum of the conformal radius of $\mathbb{D}_a$ equals $0.7$
(see (\ref{s0})), and, by Lemma~\ref{lemma2}, the maximum of the
conformal radius of $\mathbb{L}_a$, and of the symmetric to it
surface $\mathbb{L}_a^*$, is greater than $r_0=0.695356\ldots$, we
conclude that $\|B\|_{\mathbb{B}}\ge r_0$.

To prove the last part
of Theorem~\ref{Lower}, we consider $B(z)=z^n$. Then it is easy to
check that $\|z^n\|_{\mathbb{B}} \to 2/e$ as $n \to \infty$. This
fact concludes the proof of Theorem \ref{Lower}.

\medskip

\noindent {\bf Proof of Theorem~\ref{General}.} Let $f$ be
holomorphic in $\mathbb{D}$ function, $f'(0)=1$ and $B$ be a
Blaschke product of order $n$; we can assume that $n>1$. Then the
covering $f:\mathbb{D}\to f(\mathbb{D})$ defines a Riemann surface
$R$. By a well-known result (see \cite[p.364]{Ahl}), $R$ contains
a one-sheeted disk $K$ of radius $r=\sqrt{3}/4$ centered at some
point $a\in \mathbb{C}$. Therefore, there is a simply-connected
domain $G\subset \mathbb{D}$ such that $f$ maps $G$ onto $K$.
Denote by $A$ the set of all critical values of $B$. Now we
consider a sequence of disks $K_m$ centered at $a$ of radii $r_m$
such that $r_m< r$ and $r_m\to r$, $m\to \infty$. We can choose
$r_m$ such that the preimages $\gamma_m$ of the boundary circles
$\partial K_m$ under the mapping $f|G$ are disjoint with $A$. It
is evident that every $\gamma_m$ is a closed Jordan curve lying in
$G$. Now we fix $m$ and consider the set $B^{-1}(\gamma_m)$. Since
$\gamma_m$ is disjoint with $A$, the set $B^{-1}(\gamma_m)$
consists of a finite set of disjoint closed Jordan curves. We
claim that the interiors of these curves are also disjoint.
Indeed, if one of the curves, say, $\alpha$, is in the interior of
another one, $\beta$, then it can not be connected with the
boundary of $\mathbb{D}$ without intersecting $\beta$. On the
other hand, every point of $\gamma_m$ can be connected with
$\partial \mathbb{D}$ by a curve $\omega$ which does not intersect
$\gamma_m$ at other points and does not pass through points of the
set $A$. There is a unique lift of $\omega$ from some appropriate
point of $\alpha$ on $\mathbb{D}$, with respect to the covering
map $B:\mathbb{D}\to\mathbb{D}$, and the lift does not intersect
$\beta$, since $\omega$ has no common points with $\gamma$, except
for the initial one. This prove that the interiors are disjoint.
Now consider one of such curves, $\alpha$. Denote by $H$ the
interior of $\alpha$. Then $g:H\to K_m$, where $g=f\circ B$,
defines a Riemann surface which is a finite-sheeted ramified
covering of $K_m$. It is easy to see that $h(z)=(g(z)-a)/r_m$ is a
Blaschke product. Applying Theorem~\ref{Lower} to the mapping $h$,
we obtain
$$
\|g\|_{\mathbb{B}}=r_m\|h\|_{\mathbb{B}} \ge r_m\cdot r_0, \quad
m\ge 1.
$$
Taking $m\to \infty$ we obtain the desired inequality.

If $f$ is a convex univalent function, then according to \cite{Reich} the value $r=\sqrt{3}/4$ can be replaced by $\pi/4$.
Theorem~\ref{General} is proved.

\medskip

\noindent
{\bf Proof of Theorem \ref{infi}.}
Let
$$
B(z) = \prod_{j\ge 1}\frac{|z_j|}{z_j}\frac{z_j-z}{1-\bar z_j z}.
$$
Then it is easy to see that
$$
B'(z) = B(z)\sum_{j\ge 1} \frac{1-|z_j|^2}{(z-z_j)(1-\bar z_j z)}, \quad z\in \mathbb{D},
$$
and, in particular,
\begin{equation}\label{pro}
|B'(\zeta)| = \sum_{j\ge 1} \frac{1-|z_j|^2}{|\zeta -z_j|^2}.
\end{equation}

For some $d \in (0,1)$ to be chosen later, we set
$$
z_0 = \Big( 1 - \frac{d\delta}{|B'(\zeta)|} \Big) \zeta.
$$
By (\ref{onec}), $|\zeta-z_j| = |1-\bar z_j \zeta| \ge \delta |B'(\zeta)|^{-1}$ whence $z_0\in\mathbb{D}$
and the following inequalities are valid:
\begin{equation}\label{extra1}
|z-z_j| \ge |\zeta-z_j| - |\zeta-z_0|\ge (1- d) |\zeta-z_j|
\end{equation}
and
\begin{equation}\label{extra2}
|1- \overline{z_j}z| \ge |1- \overline{z_j}\zeta| - |\zeta-z_0| \ge (1- d) |\zeta - z_j|
\end{equation}
for any $z\in [z_0, \zeta]$. Therefore, $|B'(z)| \le  (1-d)^{-2} |B'(\zeta)|$, $z \in [z_0, \zeta]$. Hence,
$$
|B(z_0)-B(\zeta)| \leq \frac{|B'(\zeta)|}{(1-d)^2} \cdot |z_0-\zeta|=\frac{d\delta}{(1-d)^2}.
$$
It follows that
\begin{equation}
\label{lowerB2}
|B(z_0)| \ge 1- \frac{d \delta}{(1-d)^2}.
\end{equation}

Now we can obtain a lower estimate for $|B'(z_0)|$. We have
$$
\bigg| \frac{B'(z_0)}{B(z_0)} - \frac{B'(\zeta)}{B(\zeta)}\bigg|  =
$$
$$
 = |z_0-\zeta|\cdot \bigg| \sum_{j\ge 1}
\Big(\frac{1-|z_j|^2}{(z_0-z_j)(1-\bar z_j z_0)(\zeta-z_j)} - \frac{\bar z_j(1-|z_j|^2)}{(1-\bar z_j z_0)(\zeta-z_j)(1-\bar z_j \zeta)} \Big) \bigg|.
$$
Using (\ref{onec}), (\ref{extra1}) and (\ref{extra2}), we get
$$
\frac{1}{|(z_0-z_j)(1-\bar z_j z_0)(\zeta-z_j)|}\le \frac{1}{(1-d)^2|\zeta-z_j|^3}\le
\frac{|B'(\zeta)|}{\delta (1-d)^2 |\zeta - z_j|^2}.
$$
Estimating analogously the second term and summing up, we get
$$
\bigg| \frac{B'(z_0)}{B(z_0)} - \frac{B'(\zeta)}{B(\zeta)}\bigg|
\le
\frac{2|B'(\zeta)|^2}{\delta (1-d)^2}\cdot|z_0-\zeta| = \frac{2d |B'(\zeta)|}{(1-d)^2}.
$$
Therefore,
$$
\bigg| \frac{B'(z_0)}{B(z_0)} \bigg| \ge \bigg(1 - \frac{2d}{(1-d)^2} \bigg) |B'(\zeta)|.
$$
This estimate together with inequality (\ref{lowerB2}) yield
$$
(1-|z_0|^2) |B'(z_0)|  =(1-|z_0|^2) |B(z_0)| \cdot \bigg| \frac{B'(z_0)}{B(z_0)} \bigg| \ge
$$
$$
 \ge d\delta  \bigg( 1- \frac{d\delta}{(1-d)^2} \bigg)\bigg(1- \frac{2d}{(1-d)^2} \bigg).
$$
Recall that $\delta \le 1$. It remains to take $d= 1/7$ to obtain the required numerical
estimate. Theorem \ref{infi} is proved.
\medskip

\begin{remark}
{\rm The constants in Theorem \ref{infi} are by no means optimal. In the case when $B$ is a finite Blaschke product they can be substantially
improved. First, we can choose $\zeta$ to be a point where $|B'|$ attains its maximum in $\overline{\mathbb{D}}$. Second,
using the fact that the Bloch norm is invariant under  M\"obius transforms we can assume that the zeros are arbitrarily
close to the boundary and so $|B'(\zeta)|$ is arbitrarily large. These improvements make it possible to obtain by this method the lower
bound about 0.361..., which is still much smaller than the one obtained in Theorem \ref{Lower} by geometric methods. }
\end{remark}

\subsection*{Acknowledgments}
The authors express their deep gratitude to the anonymous referee
for the large number of useful comments that contributed to a
significant improvement of the article.

The work of I.~R. Kayumov is supported by the Russian Science Foundation under grant 18-11-00115. The work of S.R. Nasyrov performed under the development program of Volga Region
Mathematical Center (agreement no. 075-02-2021-1393).

\bibliographystyle{crplain}

\newpage

 Anton D. Baranov \\
Saint-Petersburg State University \\
Saint-Petersburg, 199178, Russia \\
e-mail: anton.d.baranov@gmail.com

Ilgiz R. Kayumov \\
Kazan Federal University \\
Institute of Mathematics and Mechanics\\
Scientific and Educational Mathematical\\ Center of the Volga Federal District\\
 Kremlevskaya 18 \\
  420 008 Kazan, Russia\\
e-mail: ikayumov@kpfu.ru

Semen R. Nasyrov \\
Kazan Federal University \\
Institute of Mathematics and Mechanics\\
 Kremlevskaya 18 \\
  420 008 Kazan, Russia\\
e-mail: semen.nasyrov@yandex.ru
\end{document}